\input amstex
\documentstyle{amsppt}
\magnification1200
\NoBlackBoxes
\def\phi{\varphi}

\def\hbar{\overline{h}}

\def\hbar{\overline{h}}
\pageheight{9 true in}
\pagewidth{6.5 true in}

\topmatter
\title
Negative values of truncations to $L(1,\chi)$
\endtitle
\author Andrew Granville and K. Soundararajan
\endauthor
\rightheadtext{Negative values}
\address{D{\'e}partment  de
Math{\'e}matiques et Statistique, Universit{\'e}
de Montr{\'e}al, CP 6128 succ
Centre-Ville, Montr{\'e}al, QC  H3C 3J7, Canada}\endaddress
\email{andrew{\@}dms.umontreal.ca}
\endemail
\address{Department of Mathematics, University of Michigan, Ann Arbor,
Michigan 48109, USA} \endaddress \email{ksound{\@}umich.edu} \endemail
\thanks{Le premier auteur est partiellement soutenu par une bourse
du Conseil  de recherches en
sciences naturelles et en g\' enie du Canada. The
second  author is partially supported by the National Science Foundation and the
American Institute of Mathematics (AIM).}
\endthanks
\endtopmatter
\def\lam{\lambda}

\document

\head 1. Introduction \endhead

\noindent Dirichlet's celebrated class number formula established that
$L(1,\chi)$ is positive for primitive, quadratic Dirichlet characters $\chi$.
One might attempt to prove this positivity by trying to establish that the partial sums
$\sum_{n\le x} \chi(n)/n$ are all non-negative.
However, such truncated sums can get negative, a feature which we will
explore in this note. 

By quadratic reciprocity we may find an arithmetic progression $\pmod {4 \prod_{p\le x}
p}$ such that any prime $q$ lying in this progression satisfies
$\fracwithdelims(){p}{q}=-1$ for each $p\le x$.  Such primes $q$ exist
by Dirichlet's theorem on primes in arithmetic progressions, and for such $q$
we have $\sum_{n\le x}
\fracwithdelims(){n}{q} /n  = \sum_{n\le x} \lam(n)/n$ where $\lam(n)= (-1)^{\Omega(n)}$
is the Liouville function.  Tur{\' a}n [6] suggested that $\sum_{n\le x}
\lam(n)/n$ may be always positive, noting that this would imply the truth
of the Riemann Hypothesis (and previously P{\' o}lya had conjectured that the related
$\sum_{n\le x} \lam(n)$ is non-positive for all $x\ge 2$, which also implies the Riemann Hypothesis).  In [4] Haselgrove showed that both the
Tur{\' a}n and P{\' o}lya conjectures are false;  therefore we know  that truncations to $L(1,\chi)$ may get negative.

Let ${\Cal F}$ denote the set of all
completely multiplicative functions $f(\cdot)$ with $-1\le f(n)\le 1$ for all
positive integers $n$, let ${\Cal F}_1$ be those for which each $f(n)= \pm 1$, and ${\Cal F}_0$ be those for which each $f(n)=0$ or $\pm 1$.
Given any $x$ and any  $f\in {\Cal F}_0$ 
we may  find a primitive quadratic character $\chi$ with $\chi(n) = f(n)$ for all
$n\le x$ (again, by using quadratic reciprocity and
Dirichlet's theorem on primes in arithmetic progressions) so that, for any $x\ge 1$, 
$$
\min\Sb\chi \ \text{a quadratic}\\ \text{character}\endSb  \ \sum_{n\le x} \frac{\chi(n)}{n} =\  \delta_0(x):=\min_{f\in {\Cal F}_0} \sum_{n\le x} \frac{f(n)}{n} .
$$
Moreover, since ${\Cal F}_1 \subset {\Cal F}_0 \subset {\Cal F}$ we have that
$$
\delta(x): = \min_{f\in {\Cal F}} \sum_{n\le x}
\frac{f(n)}{n} \quad \leq
\delta_0(x)
\leq
\quad
\delta_1(x):= \min_{f\in {\Cal F}_1} \sum_{n\le x} \frac{f(n)}{n}.
$$
We expect that $\delta(x)\sim \delta_1(x)$ and even, perhaps, that $\delta(x)= \delta_1(x)$ for sufficiently large $x$.

Trivially $\delta(x) \ge -\sum_{n\le x} 1/n = -(\log x +\gamma +O(1/x))$.  Less
trivially $\delta(x) \ge -1$, as may be shown by considering
the non-negative multiplicative function $g(n) =\sum_{d|n} f(d)$ and noting
that
$$
0\le \sum_{n\le x} g(n) = \sum_{d\le x} f(d) \Big[\frac{x}{d}\Big]
\le \sum_{d\le x} \Big(x\frac{f(d)}{d} + 1\Big).
$$
We will show that $\delta(x)\leq \delta_1(x)<0$ for all large values of
$x$, and that $\delta(x) \to 0$ as $x\to \infty$.

\proclaim{Theorem 1}  For all large $x$ and all $f\in {\Cal F}$ we have
$$
\sum_{n\le x} \frac{f(n)}{n} \ge - \frac{1}{(\log \log x)^{\frac{3}{5}}}.
$$
Further, there exists a constant $c>0$ such that   for all large $x$ there
exists a function $f (=f_x) \in {\Cal F}_1$ such
that
$$
\sum_{n\le x} \frac{f(n)}{n} \le - \frac{c}{\log x}.
$$
In other words, for all large $x$,
$$
-\frac{1}{(\log \log x)^{\frac{3}{5}}} \le \delta(x) \le \delta_0(x) \le \delta_1(x) \le
-\frac{c}{\log x}.
$$
\endproclaim

It would be interesting to determine more precisely the asymptotic
nature of $\delta(x), \delta_0(x)$ and $\delta_1(x)$, and to understand the nature of the
optimal functions.

Instead of completely multiplicative functions we may consider the
larger class ${\Cal F}^*$ of multiplicative functions, and analogously
define $\delta^*(x):= \min_{f\in {\Cal F}^*} \sum_{n\le x} f(n)/n$.

\proclaim{Theorem 2}  We have
$$
\delta^*(x) = \Big( 1 -2\log (1+\sqrt{e}) +4\int_1^{\sqrt{e}}
\frac{\log t}{t+1} dt \Big)\log 2 +o(1) = -0.4553 \ldots+o(1).
$$
If $f^*\in {\Cal F}^*$ and $x$ is large then
$$
\sum_{n\le x} \frac{f^*(n)}{n} \ge -\frac{1}{(\log \log x)^{\frac 35}},
$$
unless
$$
\sum_{k=1}^{\infty} \frac{1+f^*(2^k)}{2^k} \ll (\log x)^{-\frac 1{20}}.
$$
Finally
$$
\sum_{n\le x} \frac{f^*(n)}{n}=\delta^*(x) +o(1)
$$
if and only if
$$
\Big(\sum_{k=1}^{\infty} \frac{1+f^*(2^k)}{2^k}\Big)\log x +
\sum_{3 \le p \le x^{1/(1+\sqrt{e})}}
\sum_{k=1}^{\infty} \frac{1-f^*(p^k)}{p^k}  +
\sum_{x^{1/(1+\sqrt{e})} \le p \le x}
\frac{1+f^*(p)}{p} = o(1).
$$
\endproclaim


\head 2. Constructing negative values \endhead

\noindent Recall Haselgrove's result [4]: there exists an integer $N$ such that
$\sum_{n\le N} \lam(n)/n = -\delta$ with $\delta >0$, where
$\lambda\in {\Cal F}_1$ with $\lam(p)=-1$ for all primes $p$.
Let $x> N^2$ be large
and consider the function $f =f_x \in {\Cal F}_1$ defined by $f(p)=1$ if
$x/(N+1) < p \le x/N$ and $f(p)=-1$ for all other $p$.
If $n\le x$ then
we see that $f(n) = \lam(n)$ unless $n=p\ell$ for a (unique) prime
$p\in (x/(N+1),x/N]$ in which case $f(n) = \lam(\ell) = \lam(n) +2 \lam(\ell)$.
Thus
$$
\align
\sum_{n\le x} \frac{f(n)}{n} &=
\sum_{n\le x} \frac{\lam(n)}{n} + 2 \sum_{ x/(N+1) <p \le x/N} \frac{1}{p}
\sum_{\ell \le x/p} \frac{\lam(\ell)}{\ell}\\
&=\sum_{n\le x} \frac{\lam(n)}{n} - 2\delta \sum_{x/(N+1) < p\le x/N}
\frac{1}{p}. \tag{2.1}\\
\endalign
$$

A standard argument, as in the proof of the prime number theorem, shows
that
$$
\sum_{n\le x} \frac{\lam(n)}{n} = \frac{1}{2\pi i} \int_{2-i\infty}^{2+i\infty}
\frac{\zeta(2s+2)}{\zeta(s+1)} \frac{x^s}{s} ds \ll \exp(-c\sqrt{\log x}),
$$
for some $c>0$.  Further the prime number theorem readily gives that
$$
\sum_{x/(N+1) < p \le x/N} \frac 1p \sim \log \Big(\frac{\log (x/N)}{\log
(x/(N+1))} \Big)
\asymp \frac{1}{N \log x}.
$$
Inserting these estimates in (2.1) we obtain that $\delta(x) \le -c/\log x$
for large $x$ (here $c\asymp \delta/N$), as claimed in Theorem 1.

\head 3. The lower bound for $\delta(x)$ \endhead

\proclaim{Proposition 3.1} Let $f$ be a completely multiplicative function with
$-1\le f(n) \le 1$ for all $n$, and set $g(n) = \sum_{d|n} f(d)$ so that $g$ is
a
non-negative multiplicative function.  Then
$$
\sum_{n\le x} \frac{f(n)}{n} = \frac{1}{x} \sum_{n\le x} g(n)  + (1-\gamma)
\frac{1}{x}
\sum_{n\le x} f(n) + O\Big( \frac{1}{(\log x)^{\frac 15}}\Big).
$$
\endproclaim

\demo{Proof}   Define $F(t)= \frac{1}{t} \sum_{n\le t} f(n)$.  We will make use
of the fact that $F(t)$ varies slowly with $t$.  From [2, Corollary 3] we find
that
if $1\le w \le x/10$ then
$$
\Big| |F(x)| - |F(x/w)| \Big| \ll \Big(\frac{\log 2w}{\log
x}\Big)^{1-\frac{2}{\pi}}
\log \Big( \frac{\log x}{\log 2w} \Big) + \frac{\log \log x}{(\log
x)^{2-\sqrt{3}}}. \tag{3.1}
$$
We may easily deduce that
$$
\Big| F(x) -F(x/w) \Big| \ll \Big(\frac{\log 2w}{\log x}\Big)^{1-\frac{2}{\pi}}
\log \Big( \frac{\log x}{\log 2w} \Big) + \frac{\log \log x}{(\log
x)^{2-\sqrt{3}}}
\ll \Big(\frac{\log 2w}{\log x}\Big)^{\frac 14}.  \tag{3.2}
$$
Indeed, if $F(x)$ and $F(x/w)$ are of the same sign then (3.2) follows
at once from (3.1).  If $F(x)$ and $F(x/w)$ are of opposite signs then
we may find $1\le v\le w$ with $|\sum_{n\le x/v} f(n) |\le 1$ and then
using (3.1) first with $F(x)$ and $F(x/v)$, and second with $F(x/v)$ and
$F(x/w)$
we obtain (3.2).

We now turn to the proof of the Proposition.   We start with
$$
\sum_{n\le x} g(n) = \sum_{d\le x} f(d) \Big[ \frac{x}{d}\Big]
= x\sum_{d\le x} \frac{f(d)}{d} - \sum_{d\le x} f(d) \Big\{\frac{x}{d} \Big\}.
\tag{3.3}
$$
Now
$$
\align
\sum_{d\le x} f(d) \Big\{\frac{x}{d} \Big\} &= \sum_{j\le x } \sum_{x/(j+1)
<d\le x/j}
f(d) \Big( \frac{x}{d}-j\Big) \\
&= \sum_{j \le \log x} \int_{x/(j+1)}^{x/j} \frac{x}{t^2}
\sum_{x/(j+1)<d\le t} f(d) dt + O\Big(\frac{x}{\log x}\Big).\\
\endalign
$$
From (3.2) we see that if $j\le \log x$, and $x/(j+1)< t\le x/j$ then
$$
\sum_{x/(j+1) < d\le t} f(d) = \Big(t- \frac{x}{(j+1)} \Big) \frac{1}{x}
\sum_{n\le x}
f(n) + O\Big( \frac{x\log(j+1)}{j(\log x)^{\frac 14}}\Big).
$$
Using this above we conclude that
$$
\sum_{d\le x} f(d) \Big\{\frac{x}{d}\Big\}
= \Big(\sum_{n\le x} f(n) \Big)\sum_{j\le \log x}\Big( \log
\Big(\frac{j+1}{j}\Big) - \frac{1}{j+1}\Big)
+O\Big( \frac{x (\log \log x)^2}{(\log x)^{\frac 14}}\Big). \tag{3.4}
$$
Since $\sum_{j\le J} (\log(1+1/j)-1/(j+1)) = \log(J+1) -\sum_{j\le J+1} 1/j +1
= 1-\gamma +O(1/J)$, when we insert (3.4) into (3.3) we obtain the Proposition.

\enddemo

Set $u=\sum_{p\le x} (1-f(p))/p$.   By Theorem 2 of A. Hildebrand [5]
(with $f$ there being our function $g$, $K=2$, $K_2=1.1$, and $z=2$)
we obtain that
$$
\align
\frac{1}{x} \sum_{n\le x} g(n) &\gg \prod_{p\le x} \Big( 1-\frac 1p\Big) \Big(1+
\frac{g(p)}{p} + \frac{g(p^2)}{p^2} +\ldots \Big)
\sigma_{-} \Big( \exp\Big(\sum_{p\le x} \frac{\max(0,1-g(p))}{p} \Big)\Big)
\\
&\hskip 1 in + O(\exp(-(\log x)^{\beta})), \\
\endalign
$$
where $\beta$ is some positive constant and $\sigma_{-}(\xi) =\xi \rho(\xi)$
with $\rho$ being the Dickman function{\footnote {The Dickman function is defined as
$\rho(u)=1$ for $u\leq 1$, and $\rho(u)=(1/u)\int_{u-1}^u \rho(t) dt$
for $u\geq 1$.}}.
Since $\max(0,1-g(p)) \le (1-f(p))/2$
we deduce that
$$
\align
\frac{1}{x} \sum_{n\le x} g(n) &\gg (e^{-u} \log x) (e^{u/2} \rho(e^{u/2}))
+ O(\exp(-(\log x)^{\beta})) \\
&\gg e^{-ue^{u/2}} (\log x) + O(\exp(-(\log x)^{\beta})), \tag{3.5}
\\
\endalign
$$
since $\rho(\xi) = \xi^{-\xi +o(\xi)}$.

On the other hand, a special case of the main result in [3] implies  that
$$
\frac{1}{x }\Big| \sum_{n\le x} f(n)\Big|
\ll   e^{-\kappa u}, \tag{3.6}
$$
where $\kappa = 0.32867\ldots$.
Combining Proposition 3.1 with (3.5) and (3.6) we immediately
get that $\delta(x) \ge  -c/(\log \log x)^{\xi}$ for any
$\xi < 2\kappa$.  This completes the proof of Theorem 1.

\remark{Remark} The bound (3.5) is attained only in certain very special cases,
that is when there are very few  primes $p>x^{e^{-u}}$ for which $f(p)=1+o(1)$.
In this case one can get a far stronger bound than (3.6). Since the first part
of Theorem 1 depends on an interaction between these two bounds, this suggests
that one might be able to improve Theorem 1 significantly by determining how
(3.5) and (3.6) depend upon one another.
\endremark

\head 4. Proof of Theorem 2 \endhead

\noindent Given $f^* \in {\Cal F}^*$ we associate
a completely multiplicative function $f\in {\Cal F}$
by setting $f(p)=f^*(p)$.  We write $f^*(n)
=\sum_{d|n} h(d) f(n/d)$ where $h$ is
the multiplicative function given by
$h(p^k) = f^*(p^k)-f(p)f^*(p^{k-1})$
for $k\ge 1$.  Now,
$$
\align
\sum_{n\le x}\frac{f^*(n)}{n}
&=\sum_{d\le x} \frac{h(d)}{d} \sum_{m\le x/d} \frac{f(m)}{m}
\\
&= \sum_{d\le (\log x)^{6}} \frac{h(d)}{d} \sum_{m\le x/d} \frac{f(m)}{m}
+ O\Big( \log x \sum_{d> (\log x)^{6} }\frac{|h(d)|}{d} \Big). \tag{4.1}
\\
\endalign
$$
Since $h(p)=0$ and $|h(p^k)| \le 2$ for $k\ge 2$ we
see that
$$
\sum_{d> (\log x)^6} \frac{|h(d)|}{d}
\le (\log x)^{-2} \sum_{d \ge 1} \frac{|h(d)|}{d^{\frac 23}} \ll (\log x)^{-2}.
\tag{4.2}
$$
Further, for $d\le (\log x)^6$, we have (writing $F(t)=\frac 1t \sum_{n\le t}
f(n)$ as
in section 3)
$$
\sum_{x/d \le n\le x} \frac{f(n)}{n}
= F(x) - F(x/d) + \int_{x/d}^{x} \frac{F(t)}{t} dt
= \frac{\log d}{x} \sum_{n\le x} f(n) +O\Big(\frac{1}{(\log x)^{\frac 15}}\Big),
$$
using (3.2). Using the above in (4.1) we
deduce that
$$
\sum_{n\le x} \frac{f^*(n)}{n}
= \Big( \sum_{n\le x} \frac{f(n)}{n}
\Big) \sum_{d\le (\log x)^6} \frac{h(d)}{d}
- \frac{1}{x} \sum_{n\le x} f(n) \sum_{d\le (\log x)^6} \frac{h(d)\log d}{d}
+ O\Big(\frac{1}{(\log x)^{\frac 15}}\Big).
$$
Arguing as in (4.2) we may extend the sums over $d$ above to all $d$,
incurring a negligible error.  Thus we conclude that
$$
\sum_{n\le x} \frac{f^*(n)}{n}
= H_0 \sum_{n\le x} \frac{f(n)}{n}
+ H_1 \frac{1}{x} \sum_{n\le x} f(n) + O\Big( \frac{1}{(\log x)^{\frac
15}}\Big),
$$
with
$$
H_0 =\sum_{d=1}^{\infty} \frac{h(d)}{d},
\qquad \text{and} \qquad H_1
= -\sum_{d=1}^{\infty} \frac{h(d)\log d}{d}.
$$
Note that $H_0 = \prod_{p} (1+h(p)/p+h(p^2)/p^2 +\ldots) \ge 0$, and that $H_0,
|H_1|\ll 1$.

We now use Proposition 3.1, keeping the notation there.  We deduce
that
$$
\sum_{n\le x} \frac{f^*(n)}{n}
= H_0 \frac{1}{x} \sum_{n\le x} g(n)
+ \Big( (1-\gamma)H_0 + H_1 \Big) \frac{1}{x} \sum_{n\le x} f(n)
+ O\Big( \frac{1}{(\log x)^{\frac 15}} \Big).  \tag{4.3}
$$
If $H_0 \ge (\log x)^{-\frac 1{20}}$ then we may argue as in section 3,
using (3.5) and (3.6).  In that case, we see that $\sum_{n\le x} f^*(n)/n
\ge -1/(\log \log x)^{\frac 35}$.  Henceforth we suppose that
$H_0 \le (\log x)^{-\frac 1{20}}$. Since
$$
H_0 \asymp 1+ \frac{h(2)}{2}+ \frac{h(2^2)}{2^2} +\ldots
\asymp 1+ \frac{f^{*}(2)}{2} +\frac{f^*(2^2)}{2^2} + \ldots,
$$
we deduce that (note $h(2)=0$)
$$
\sum_{k=2}^{\infty} \frac{2+h(2^k)}{2^k}
\asymp \sum_{k=1}^{\infty} \frac{1+f^*(2^k)}{2^k} \ll (\log x)^{-\frac{1}{20}}. \tag{4.4}
$$
This proves the middle assertion of Theorem 2.





Writing $d =2^k \ell$ with $\ell$ odd,
$$
\align
H_1 &= -\sum_{ \ell \text{ odd}} \frac{h(\ell)}{\ell}
\sum_{k=0}^{\infty} \frac{h(2^k)}{2^k} (k\log 2 +\log \ell)
\\
&= - \log 2 \Big(\sum_{k=1}^{\infty} \frac{k h(2^k)}{2^k} \Big) \sum_{\ell \text{ odd}}
\frac{h(\ell)}{\ell} + O((\log x)^{-\frac 1{20}})\\
&= 3\log 2 \prod_{p\ge 3} \Big(1+\frac{h(p)}{p} +\frac{h(p^2)}{p^2}+\ldots \Big)
+ O\Big(\frac{\log \log x}{(\log x)^{\frac{1}{20}}}\Big),\\
\endalign
$$
where we have used (4.4) and that $\sum_{k=1}^{\infty} kh(2^k)/2^k
= -3 +O(\log \log x/(\log x)^{\frac 1{20}})$.  Using these observations in (4.3)
we obtain that
$$
\align
\sum_{n\le x} \frac{f^*(n)}{n}
&= H_0 \frac 1x \sum_{n\le x} g(n) + 3\log 2 \prod_{p\ge 3}
\Big(1+\frac{h(p)}{p} +\frac {h(p^2)}{p^2} +\ldots\Big) \frac{1}{x}
\sum_{n\le x} f(n) + o(1) \\
&\ge  3\log 2 \prod_{p\ge 3}
\Big(1+\frac{h(p)}{p} +\frac {h(p^2)}{p^2} +\ldots\Big) \frac{1}{x}
\sum_{n\le x} f(n) + o(1). \tag{4.5}
\\
\endalign
$$

Let $r(\cdot)$ be the completely multiplicative function
with $r(p)=1$ for $p\leq \log x$,  and $r(p)=f(p)$ otherwise.  Then
Proposition 4.4 of [1] shows that
$$
\frac{1}{x} \sum_{n\le x} f(n) = \prod_{p\leq \log x}
\Big( 1- \frac 1p\Big) \Big( 1- \frac {f(p)}p\Big)^{-1} \frac{1}{x}
\sum_{n\le x} r(n) + O\Big( \frac 1{(\log x)^{\frac{1}{20}} } \Big).
$$
Since $f(2)=-1+O(H_0)$ we deduce from (4.5) and the above that
$$
\sum_{n\le x} \frac{f^*(n)}{n}
\ge \log 2 \prod_{p\ge 3} \Big(1-\frac 1p\Big) \Big( 1+
\frac{f^*(p)}{p} + \frac{f^*(p^2)}{p^2}+\ldots \Big)
\frac{1}{x} \sum_{n\le x} r(n) + o(1). \tag{4.6}
$$

One of the main results of [1] (see Corollary 1 there) shows that
$$
\frac 1x \sum_{n\le x} r(n) \ge 1 -2\log (1+\sqrt{e}) + 4\int_1^{\sqrt{e}}
\frac{\log t}{t+1} dt + o(1)
= -0.656999\ldots +o(1), \tag{4.7a}
$$
and that equality here holds if and only if
$$
\sum_{p\le x^{1/(1+\sqrt{e})}} \frac{1-r(p)}{p} +
\sum_{x^{1/(1+\sqrt{e})} \le p\le x} \frac{1+r(p)}{p}
=o(1). \tag{4.7b}
$$
Since the product in (4.6) lies between $0$ and $1$ we conclude
that
$$
\sum_{n\le x} \frac{f^*(n)}{n}
\ge \Big(1 -2\log (1+\sqrt{e}) + 4\int_1^{\sqrt{e}}
\frac{\log t}{t+1} dt\Big) \log 2 +o(1), \tag{4.8}
$$
and for equality to be possible here we must have (4.7b), and
in addition that the product in (4.6) is $1+o(1)$.  These conditions
may be written as
$$
\sum_{3\le p \le x^{1/(1+\sqrt{e})} } \sum_{k= 1}^{\infty} \frac{1-f^*(p^k)}{p^k}
+ \sum_{x^{1/(1+\sqrt{e})} \le p\le x} \frac{1-f^*(p)}{p} = o(1).
$$
If the above condition holds then, by (3.5), $\sum_{n\le x} g(n)
\gg x \log x$ and so for equality to hold in (4.5) we must have
$H_0 = o(1/\log x)$.  Thus equality in (4.8) is
only possible if
$$
\Big(\sum_{k=1}^{\infty} \frac{1+f^*(2^k)}{2^k} \Big)
\log x + \sum_{3\le p \le x^{1/(1+\sqrt{e})} } \sum_{k= 1}^{\infty} \frac{1-f^*(p^k)}{p^k}
+ \sum_{x^{1/(1+\sqrt{e})} \le p\le x} \frac{1-f^*(p)}{p} = o(1).
$$
Conversely, if the above is true then equality holds in (4.5), (4.6), and
(4.7a) giving equality in (4.8).  This proves Theorem 2.

\enddocument